\documentclass[12pt]{amsart}

\usepackage{graphicx}
\usepackage{amssymb}
\usepackage{amsfonts}

\usepackage{times,amsmath,amsthm,amssymb,latexsym,amscd}

\newtheorem{Th}{Theorem}[section]

\newcommand{\btbt}{\left( \begin{array}{cc}}
\newcommand{\etbt}{ \end{array}\right)}
\newcommand{\bthbth}{\left( \begin{array}{ccc}}
\newcommand{\ethbth}{ \end{array}\right)}
\newtheorem{Lm}{Lemma}[section]

\newcommand{\bcol}{\left(\begin{array}{c}}
\newcommand{\ecol}{\end{array}\right)}

\title{An elementary proof that random Fibonacci sequences grow exponentially.}
\author{Eran Makover and Jeffrey McGowan}\thanks{We would like to thank the referee for very helpful comments.}
\address{Dept. of Mathematics\\Central Connecticut State University\\New Britain, CT}
\begin{document}
\begin{abstract} We consider random Fibonacci sequences given by $x_{n+1}=\pm \beta x_{n}+x_{n-1}$.  Viswanath (\cite{viswanath}), following Furstenberg (\cite{furst}) showed that when $\beta = 1$, $\lim_{n\to \infty}|x_{n}|^{1/n}=1.13\ldots$ , but his proof involves the use of floating point computer calculations.  We give a completely elementary proof that  $1.23375 \ge (E(|x_{n}|))^{1/n} \ge 1.12095$ where $E(|x_{n}|)$ is the expected value for the absolute value of the $n$th term in a random Fibonacci sequence.  We compute this expected value using recurrence relations which bound the sum of all possible $n$th terms for such sequences. In addition, we give upper ands lower bounds for the second moment of the $|x_{n}|$.   Finally, we consider the conjecture of Embree and Trefethen (\cite{et}), derived using computational calculations, that for values of $\beta < 0.702585$ such sequences decay.  We show that as $\beta$ decreases, the critical value where growth can change to decay is in fact $\frac{1}{\sqrt{2}}$.\end{abstract}
\maketitle

Consider the {\it random Fibonacci sequence} generated by the recursive rule $x_{n+1}=\pm x_{n}+x_{n-1}$.   Clearly, if one considers the expected value $E(x_{n})$, this will be equal to $x_{0}$, but if one considers instead $E(|x_{n}|)$ this is not the case.  Here,  $E(|x_{n}|)$ is the expected value for the absolute value of the $n$th term in a random Fibonacci sequence.  While it might seem possible that somehow plusses and minuses cancel, and the average random Fibonacci sequence remains bounded, Viswanath (\cite{viswanath}) showed that in fact almost every sequence grows exponentially with mean value $\approx {(1.13\ldots )^{n}}$.  For more general random sequences  $x_{n+1}=\pm \beta x_{n}+x_{n-1}$, Embree and Trefethen gave numerical evidence for the dividing line between sequences which grow and those which decay (\cite{et}) - this case has also been investigated by Sire and Krapivsky (\cite{sk}).

Unfortunately, Viswanath's proof is quite complex, and involves the use of floating point computer calculations.  In this note, we give a proof that is completely elementary, and while we don't obtain as precise a value for the growth rate, we do obtain bounds on the variance for $x_{n}$ as well. 

We consider the set of possible random Fibonacci sequences as a tree, with the first leaf being $x_{2}$ (which we label $n=0$), and each leaf $x_{i}$ having two branches leading to the child leaves $x_{i-1}+x_{i}$ and $x_{i-1}-x_{i}$.   Such a tree lists the $2^{n}$ such sequences of length $n$.  At any given level in the tree, we can replace negative values with their absolute values which amounts to nothing more than a reversal of the orientation of the leaves (which is addition and which subtraction).   We will develop recursive formulas for the sum of the absolute values of the entries in a row in terms of the sums of the previous three rows - the reorientation will allow us to assume that the values in our initial two rows are positive.  The expected value $E(|x_{n}|)$ can then be bounded since it is simply the sum of the $n$th row divided by $2^{n}$.

\begin{figure}[htbp] 

   \includegraphics[width=6in]{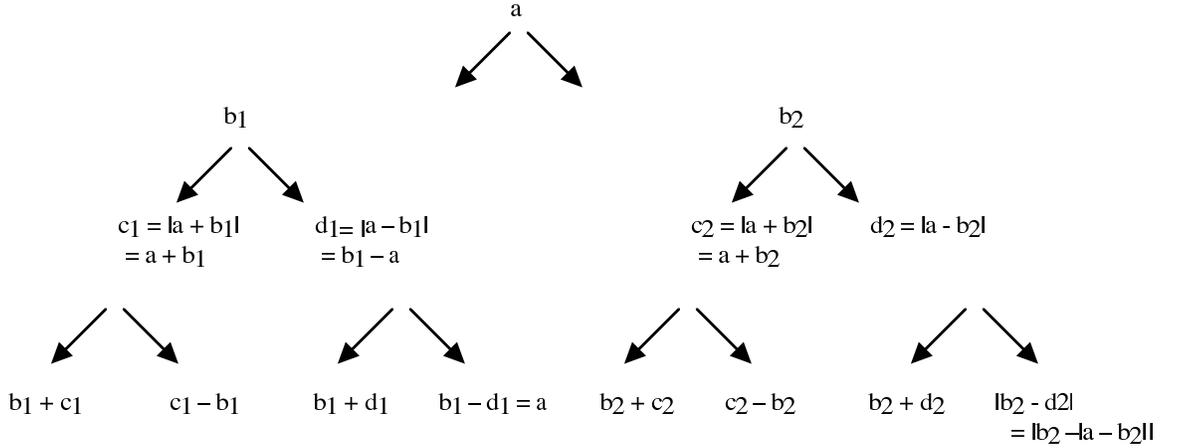} 

   \caption{The tree below $a$.  We assume $b_{1}\ge a$.}
   \label{treepic}
\end{figure}

We begin by considering the leaves below a given entry $a$ in the tree (see Figure \ref{treepic}).  $a$ has two children and there are two possibilities for each, either $b \ge a$  or $b < a$ .   It is clearly impossible that both children can be less than the parent, and we assume that $b_{1}\ge a$, but we make no assumptions yet about the value of $b_{2}$ .

We consider the sum of the last row, and bound this above and below by multiples of the sums of the previous three rows.  We note that clearly it would be possible to eliminate any $c$'s and $d$'s from Figure \ref{treepic}, but bounds obtained in this way are not as good as the bounds we give below.  The only term in the sum of the bottom row which depends on the value of $b_{2}$ is $|b_{2} - |a - b_{2}||$.

The sum of the last row is given by $$\sigma = a + b_{1}+b_{2}++2c_{1}+2c_{2}+d_{1}+d_{2}+|b_{2}-|a-b_{2}||.$$  We have the following upper and lower bounds for this sum,
\begin{Lm}\label{nicebounds}$$4a + b_{1}+b_{2}+c_{1}+c_{2}+d_{1}+d_{2}\le \sigma \le 4a + 2b_{1}+2b_{2}+c_{1}+c_{2}+d_{1}+d_{2}.$$\end{Lm}
\begin{proof}
The left hand inequality is immediate, since $c_{1}+c_{2}=2a + b_{1} + b_{2}$ and we assume that $b_{1}\ge a$.  To show the right hand inequality, we need only show that $|b_{2}-|a-b_{2}|| \le a$.  We need to consider three cases.
\begin{enumerate}
\item $b_{2}\ge a$.  Then $|b_{2}-|a-b_{2}|| = |b_{2}-(b_{2}-a)| = a$.
\item $ a > b_{2} \ge a/2$.  Then $|b_{2}-|a-b_{2}|| = |b_{2}-(a - b_{2})| = |2b_{2}-a| =2b_{2}-a <a$.
\item $b_{2} < a/2$.  Then $|b_{2}-|a-b_{2}|| = |b_{2}-(a - b_{2})| = |2b_{2}-a| =a - 2b_{2} <a$.
\end{enumerate}\end{proof}

If we denote the sum of the $n$th row by $S[n]$, then the lower bound in Lemma \ref{nicebounds} can be written as $$4S[n-3]+S[n-2]+S[n-1]$$ while the upper bound is $$4S[n-3]+2S[n-2]+S[n-1].$$  These give recurrence relations for the sums of the rows which will give upper and lower bounds for the actual sums.  If we assume that our sequence of sums consists of rational numbers, solutions to the relations can be obtained which depend on the real roots of the irreducible cubics $x^{3}-x^{2}-x-4$ and $x^{3}-x^{2}-2x-4$.  The growth rates for $S[n]$ will be given by these roots divided by 2 (to account for the doubling of the number of entries in each row.

We get  \begin{Th}\label{th1}$$1.23375\ldots\ge  (E(|x_{n}|))^{\frac{1}{n}}\ge 1.12095\ldots$$\end{Th}
   
 To determine bounds for the variance, we need to consider the sum of the {\emph squares} of a given row (the second {\it raw moment} $\mu_{2}'$), which we label $SS[n]$.  This is considerably easier.  If we consider starting again with some entry $a$, with children $b_{1}$ and $b_{2}$, then the grandchildren of $a$ will have absolute values $|a-b_{1}|$,$|a+b_{1}|$,$|a-b_{2}|$,$|a+b_{2}|$.  Squaring eliminates all the cross terms, and we just get a sum of squares of $4a^{2}+2b_{1}^{2}+2b_{2}^{2}$, which is four times the sum of the squares of the first row plus twice the sum of squares of the second.  The recurrence relation $$SS[n]=2SS[n-1]+4SS[n-2]$$ gives a growth factor for $\mu_{2}'$ of $1+\sqrt{5}$.  Now, the variance is given by $\mu_{2}=-\mu_{1}'^{2}+\mu_{2}'$ and we have upper and lower bounds on $\mu_{1}'$ using the growth factors in  Table 1.  In either case, as $n \to \infty$, the contribution from the $\mu_{1}'^{2}$ disappears.  This gives \begin{Th}$$\lim_{n \to \infty} (\mu_2|x_n|)^\frac{1}{n}=1+\sqrt{5}$$\end{Th}
 
\begin{figure}[htbp]
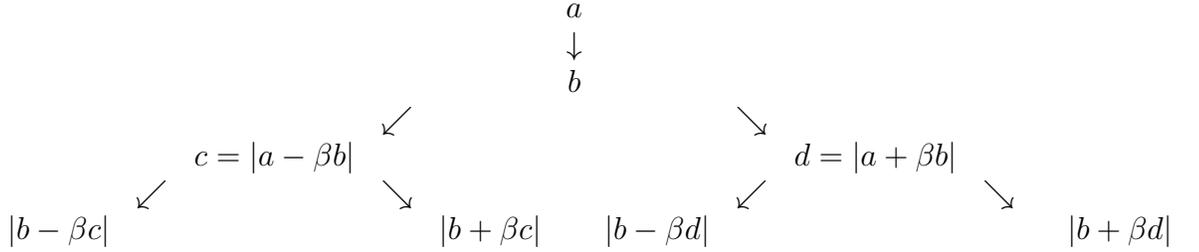
 
   \centering
$$\begin{array}{cccccccccc}  &   &   &   & a &   &   &   &   &   \\  &   &   &   & \downarrow &   &   &   &   &   \\  &   &   &   & b &   &   &   &   &   \\  &   &   & \swarrow &  & \searrow &   &   &   &   \\  &   & c =|a-\beta b|&   &  &   & d =|
a+\beta b|&   &   &   \\  & \swarrow &   & \searrow &   & \swarrow &   & \searrow &   &   \\|b-\beta c| &   &   &   & |b + \beta c|\qquad |b-\beta d| &   &   &   &   & |b +\beta d|  \\  &   &   &   &   &   &   &   &   &   \\  &   &   &   &   &   &   &   &   &   \\  &   &   &   &   &   &   &   &   &  \end{array}$$
   \caption{Half the tree below $a$ when a multiplicative constant is included.}
   \label{treepic2}
\end{figure}

 Finally, if we generalize the allowable sequences by allowing multiplication by a factor $\beta$, we get a modified version of Figure \ref{treepic}, shown in Figure \ref{treepic2}.  Here we consider only half the tree below $a$, and account for possible combinations of these halves below.
   There are two possible sums for the absolute values of the bottom row,  depending on the sign in the  absolute value on the left in the previous row, \begin{eqnarray}
\beta b \ge a & \implies & \beta (c + d) + 2b + |b + \beta a - \beta^{2} b| +  |b - \beta a - \beta^{2} b| \\
\beta b < a & \implies &\beta (c + d) + 2b + |b +\beta^{2} b - \beta a| +  |b - \beta^{2} b - \beta a|
\end{eqnarray}  There are now various possibilities for these sums, depending on whether the number in either or both absolute values is negative.  This gives six possible sums, listed with the conditions which generate them and the restrictions they impose in Table \ref{sumtable}.

\begin{table}[htdp]
  \centering \begin{tabular}{ccc}Conditions & Restrictions & Sum \\\hline &  &  \\1.  $\beta b \ge a $, $\beta a +\beta^2 b \ge b $ & $\beta^2 > \frac{1}{2}$ & $\beta(c + d) + 2b + 2\beta a$ \\$2.  \beta b \ge a $, $\beta^2 b > b + \beta a$ & $\beta > 1$ & $\beta(c + d) + (2+2\beta^2)b$ \\3.  $\beta b \ge a $, $\beta a +\beta^2 b < b $ & $\beta < 1$ & $\beta(c + d) + (4 - 2\beta^2 b$ \\4.  $\beta b < a $, $\beta a +\beta^2 b \ge b $ &  & $\beta(c + d) + (2+2\beta^2)b$ \\5.  $\beta b < a $, $\beta a> \beta^2 b + b $ & $b=0$ & $\beta(c + d) + 2\beta a$ \\6.  $\beta b < a $, $\beta a +\beta^2 b < b $ & $\beta^2 < \frac{1}{2}$ & $\beta(c + d) + 2b - 2\beta a$\end{tabular}\caption{}
\label{sumtable}
\end{table}
Note that if one combines these half sums, a number of combinations are impossible because of the conditions for such a sum to occur.  It is possible that both children of $a$ will satisfy one set of conditions, thus combinations do include some use of the same row twice. Because of the subtraction, the sixth row must be involved if the expected value is to decay exponentially, and this means that exponential decay can only occur if the condition for this is satisfied, namely if $$\beta^{2}<\frac{1}{2} \implies \beta < \frac{1}{\sqrt{2}} \approx  0.707107.$$  One would not expect decay to begin immediately as $\beta$ crosses this value, since the growth or decay rate will be a combination of various possible relations.x

\begin{figure}[htb] 
   \centering
   \includegraphics[width=6in]{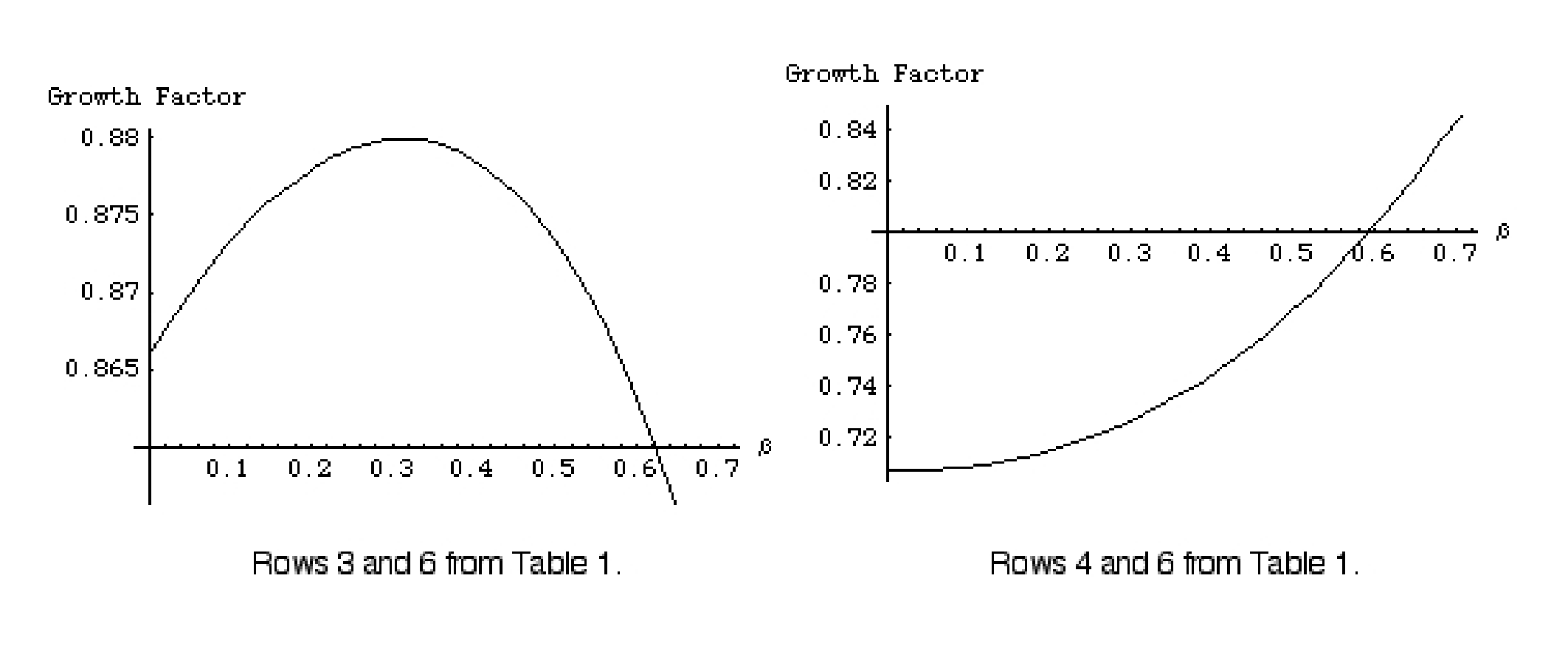} 
   \caption{Possible growth and decay as a function of $\beta$.}
   \label{growthgraphs}
\end{figure}

\begin{figure}[htb] 
   \centering
   \includegraphics[width=6in]{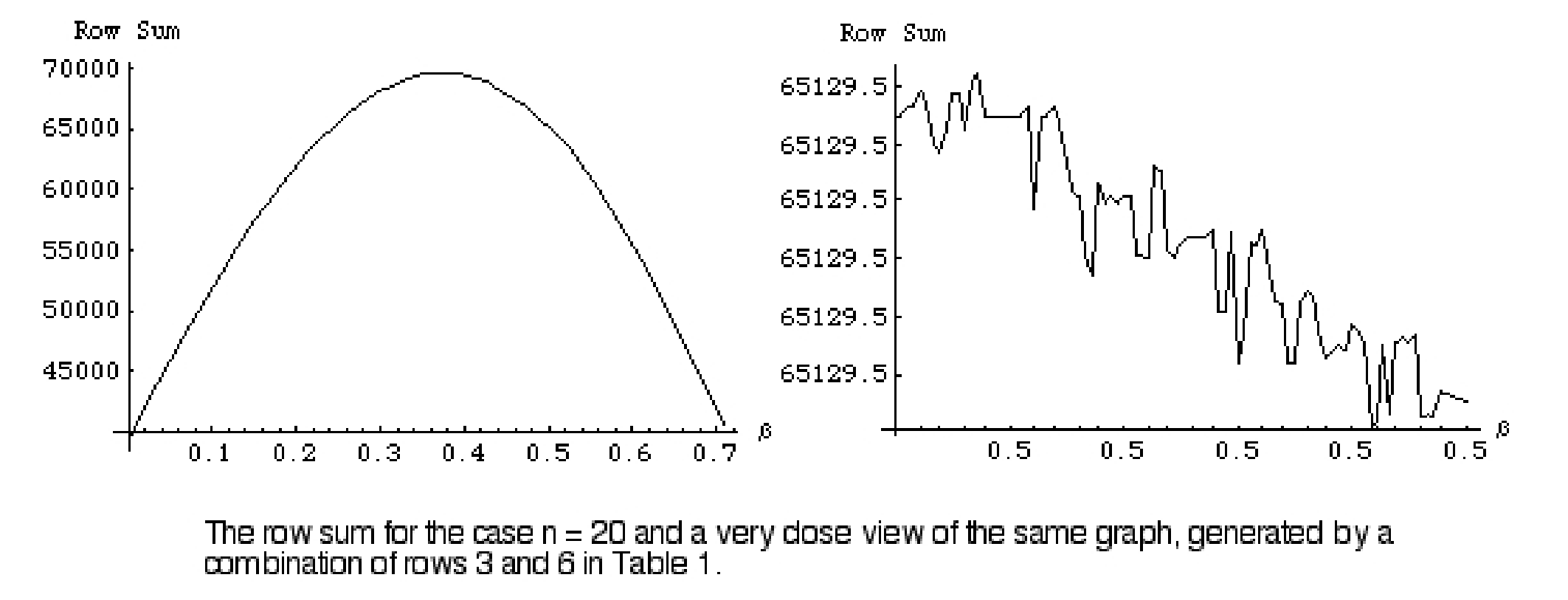} 
   \caption{Mean value of the row sum for a fixed row index $n$.}
   \label{betagraphs}
\end{figure}

We show graphs for some other possible growth factors for values of $\beta$ between 0 and 1 in Figure \ref{growthgraphs}.   If, instead, one fixes a row index $n$, it is clear that the expected value of the sum of the $n$th row does {\it not} depend smoothly on $\beta$, although one needs to ``zoom in'' quite far to see this for any reasonable large values of $n$ (Figure \ref{betagraphs}).  It seems possible that further investigation of the probability for each row to occur might allow for a deeper understanding of such scaled random Fibonacci sequences.

\bibliographystyle{amsplain} 
\bibliography{newrandomfib}
 \end{document}